\input mssymb
\footline={\hfill\number\folio\hfill\jobname\hskip5mm\number\day.\number\month.\number\year}
\def\rest{\mathord{\restriction}}
\def\ms{\medskip}

\def\ss{\smallskip}

\def\phi{\varphi}

\def\su{\subseteq}
\def\a{\alpha}
\def\b{\beta}

\def\e{\nu}
\def\d{\delta}
\def\l{\lambda}

\def\om{\omega}

\def\lng{\langle}
\def\rng{\rangle}
\def\ov{\overline}
\def\sm{\setminus}
\def\nac{{\rm nacc }\,}
\def\nacc{\nac}
\def\acc{{\rm acc}\,}

\def\cf{{\rm cf}}

\def\epsilon{\nu}
\def\otp{{\rm otp}\,}
\baselineskip=16pt

\def\endproof{\hfill $\dashv$}

\def\proclaim#1{\noindent{\bf #1}: }
\def\imply{\Rightarrow}

\def\Lev{{\rm Lev}}

\def\proof{\smallbreak\noindent{\sl Proof}: }

\magnification=1100
\def\endproof{$\dashv$}
\centerline{{\bf ${\bf \mu}$-Complete Souslin Trees on ${\bf
\mu^+}$}}
\ms
\centerline {by Menachem Kojman and Saharon
Shelah\footnote\dag{The second author thanks the Binational Science
Foundation for supporting this research. Publication no. 449}}
\centerline{The Hebrew University, Jerusalem}

\centerline{ABSTRACT}

We prove that $\mu=\mu^{<\mu}$,  $2^\mu=\mu^+$ and ``there is a non
reflecting stationary subset of $\mu^+$ composed of ordinals of
cofinality $<\mu$'' imply that there is a
$\mu$-complete Souslin tree on $\mu
^+$. 
\bigbreak

{\bf Introduction} The old problem of the existence of Souslin trees
has attracted the attention of many (see [Je] for history). While the
$\aleph_1$ case is settled, the consistency of
$GCH+SH(\aleph_2)$ is still an
open question. Gregory showed in [G]  that $GCH$ + ``there is a non
reflecting stationary set of $\om$-cofinal elements of $\om_2$''
implies the existence of an $\aleph_2$-Souslin tree. Gregory's result
showed that the consistency strength of $GCH+SH(\aleph_2)$ is at least
that of the existence of a Mahlo cardinal. Without $GCH$, the
consistency of $CH + SH(\aleph_2)$ is known from [LvSh 104]. In [ShSt
279] the equiconsistency of the existence of a weakly compact cardinal
with ``every $\aleph_2$-Aronszajn tree is special'' is shown. In [ShSt
154] it is shown that under $CH$, the consistency strength of ``there
are no $\aleph_1$-complete $\aleph_2$-Souslin trees'' is at least that of an
inaccessible cardinal.

 We show how a Souslin tree which is $\mu$-complete ($\mu$ regular)
can be constructed on a cardinal $\mu^+$ from a certain combinatorial
principle (Theorem 2 below), and then show how this principle may be
gotten from $GCH$ and a non reflecting stationary set of ordinals with
cofinality $<\mu$ in $\mu^+$ (Theorem 3 below). As a corollary
(Corollary 5 below), $GCH$ + ``there is a non reflecting stationary
set of $\om$-cofinal elements of $\om_2$'' implies the existence of an
$\aleph_1$-complete Souslin tree on $\aleph_2$.

 As mentioned in [G], 1.10(3), $CH$ and the existence of a diamond
sequence on $\{\d<\aleph_2:\cf\d=\aleph_1\}$ imply the existence of a
Souslin tree on $\aleph_2$ which is $\aleph_1$-complete. The
construction of such a tree is by induction on levels, where at a stage
of countable cofinality all branches are realized (for
$\aleph_1$-completeness), while at stages of cofinality $\aleph_1$
the diamond is consulted to realize only a part of the cofinal
branches in a way which kills all future big antichains.  The
combinatorial principle used in Theorem 2 to construct a
$\mu$-complete Souslin tree on $\mu^+$ under $GCH$ can be viewed as a
weaker substitute for a diamond sequence on $\{\d<\mu^+:\cf\d=\mu\}$:
instead of using a guess at a single guessat stage of cofinality $\mu$, we use
unboundedly many guesses, each at a level of cofinality $<\mu$. 

Unlike a diamond sequence on the stationary set of critical
cofinality, this principle makes sense also in the case of an
inaccessible cardinal (where there is no ``critical cofinality''). 
This principle is closely 
related to club guessing (see [Sh-g] and [Sh-e]), which was discovered
while the second author was trying to prove some results in Model
Theory. This principle continues the principle that
appears in [AbShSo 221], in which Souslin trees on successors of
singulars are treated. 

We learned from the referee that Gregory presented in the seventies in
a seminar at Buffalo a construction of  a countably complete
Souslin tree on $\aleph_2$ from GCH and a square, but that this was
not written. 

\proclaim {1. Notation}: (1) If $C$ is a set of ordinals, then $\acc C$
is the set of accumulation points of $C$ and $\nac C=^{df} C\sm \acc
C$. By $T_\a$ we denote the $\a$-th level of the tree $T$ and by
$T(\a)$ we denote $\bigcup_{b<\a}T_\b$.

\medbreak

\proclaim {2. Theorem} Suppose that 
\item{(a)} $\l=\mu^+=2^\mu$, $\mu=\mu^{<\mu}$;
\item{(b)} $S^*\su\{\a\in\l:\cf\a=\mu\}$ and $\ov C=\lng
c_\d:\d\in S^* \}$, $\d=\sup c_\d$,
$c_\d$ is a closed set of limit ordinals. 
\item {(c)} For every $\d\in S^*$ and
$\a\in\nac c_\d$, $P_{\d,\a}\su{\cal P}(\a)$, $|P_{\d,\a}|\le\cf\a$, and
if $\a\in S^*$, then $|P_{\d,\a}|<\mu$;
\item {(d)} For every set $A\su\l$ and club $E\su\l$, there is
a stationary $S_{A,E}\su S^*$
such that for every $\d\in S_{A,E}$, $\d=\sup\{\a\in\nac c_\d:A\cap\a\in
P_{\d,\a}\wedge \a\in E\}$;
\item{(e)} If $\d,\d^*\in S^*$, $\d\in\acc c_{\d^*}$,
then there is some $\a<\d$ such that $\lng P_{\d,\b}:\b\in \nac
c_{\d}\wedge \b>\a\rng=\lng P_{\d^*,\b}:\b\in (\nacc c_\d\cap
\d)\wedge \b>\a\rng$.
\item{(f)} for every $\gamma<\l$, $|\{\lng P_{\d,\a}:\a\in
C_\d\cap\gamma\rng:\d\in S\}|\le \mu$.

{\it Then} there is a $\mu$-complete Souslin tree on $\l$.

\proclaim{Discussion}
Condition (d) is the prediction demand. It says that for every club
$E$ and a set $A$ there is a stationary set of $\d$-s, such that for
unboundedly many non-accumulation points $\a$ of $c_\d$ two things
happen: $\a\in E$ {it and} $A\cap \a$ is guessed by $P_{\d,\a}$.


\proclaim{Proof} We assume, without loss of generality, that
for every $\d\in S^*$ and $\a\in c_\d$, $\a=\mu\a$. By induction on
$\a<\l$ we construct a tree $T(\a)$ of height $\a$ such that:
\item{(i)} The universe of $T(\a)$ is $\mu(\a+1)$, the $\b$-th level in
$T(\a)$, $T_\b$, consists of the elements $[\mu\b,\mu(\b+1))$, and
every $x\in T_\b$ for $\b<\a$ has an extension in  $T_\gamma$
for every $\gamma<\a$. Every $x\in T(\a)$ such that $\Lev(x)+1<\a$ has
at least two immediate successors in $T_\a$.
\item {(ii)} $T(\a)$ is $\mu$-complete;
\item {(iii)} For $\a<\b$, $T(\a)=T(\b)\rest |T(\a)|$

Also, we define a partial function (which, intuitively speaking,
chooses branches which help us in preserving the maximality of small
antichains that occur along the way): 
\item{(iv)} For every $x\in T(\a)$ and a
sequence $t=\lng P_{\d,\b}:\b\in\nacc c_\d\cap \a \rng$ such that
$\sup(c_\d\cap\a)<\a$ and $\Lev(x)<\max(c_\d\cap\a)$ and
$\max(c_\d\cap\a)\in\nac c_\d$, $y(x,t)$ is defined, and is an element
in the level $\max(c_\d\cap\a)$ which extends $x$ and has the property
that for every $A\in P_{\d,\max(c_\d\cap\a)}$ which is a maximal
antichain of $T({\max c_\d\cap\a})$, there is an element of $A$ below
$y(x,t)$. 
\item{(v)} If the sequence $s$ extends
the sequence $t$ and $y(x,t),y(x,s)$ exist, then $T(\a)\models y(x,t)<y(x,s)$.
\item {(vi)} For every
increasing sequence $\lng t_i:i<i^*\rng$ there is an upper bound (in
the tree order) to $\lng y(x,t_i):i<i^*\rng$.

The last demand is: 
\item{(vii)} If $\a=\d+1$, $\d\in S^*$ then
every $y\in T_\d$  satisfies that there is some
$\d^*\ge\d\in\acc c_{\d^*}$ and $x\in T(\d)$, such that $y$ is the
least upper bound (in the tree order) of $\lng y(x,t_\a):\a\in\nac
c_{\d^*}\cap\d\,\wedge\,\a_x<\a<\d\rng$ where $\a_x$ is the least in
$\nac c_{\d^*}$ such that $\a_{x}>\Lev(x)$, and $t_\a=\lng
P_{\d,\b}:\b\in \nac c_{d^*}\wedge \b\le \a\rng$.

We first show that this construction, once carried out, yields a
$\mu$-complete Souslin tree on $\l$. The completeness of $T=\cup T(\a)$
is clear from the regularity of $\l$. Suppose that $A\su \l$ is a
maximal antichain of $T$ of size $\l$. Let $E$ be the club of points
$\d<\l$ such that $T\rest\d=T(\d)$ and $A\rest \d$ is a maximal
antichain of $T(\d)$.  Pick a point $\d\in S^*$ such that
$\d=\sup\{\a\in\nacc c_\d:\a\in E\wedge A\rest \a\in P_{\d,\a}\}$. As
$|T(\d)|<\l$ there is an element $a\in A$, $\Lev(a)>\d$.  Let $y$ be
the unique such that $\Lev(y)=\d$ and $y<a$. Then by demand (vii),
there is some $\d^*\ge\d$ and $x\in T(\d)$ such that $y$ is the least
upper bound (in the tree order) of $\lng y(x, t_\a):\a\in \nac
c_{\d^*}\cap\d\wedge
\a>\Lev(x)\rng $. There is some $\a^*<\d$ such that $\lng
P_{\d,\b}:\a<\b<\d\wedge\b\in \nac c_\d\rng=\lng P_{\d^*,\b}:\a^*<\b\in \nac
c_{\d^*}\cap \d\rng$.  Pick some $\a\in \nac c_\d$ such that
$\a>\max\{\Lev (x), \a^*\}$, $\a\in E$ and $A\rest \a\in
P_{\d,\a_i}$. So $\a\in \nacc c_{\d^*}$ and $A\cap \a\in
P_{\d^*,\a}$. Then the unique $x'<y$ with $\Lev(x')=\a$ (which
equals $y(x,\lng P_{\d^*,\gamma}:\gamma\in (\nacc c_{\d^*}\cap
(\a+1))\rng)$) is above an element $a'\in A\rest\a$.  But $x'<a$
--- a contradiction to the fact that $A$ is an antichain.

Next let us show that we can carry out the construction by induction.
When $\a=\b+1$ and $\b$ is a successor or zero, add two immediate successors
to every point in the $\b$-th level. When $\b$ is limit, $\cf\b<\mu$,
add an element above every infinite branch. This addition amounts to
the total of $\mu^{<\mu}=\mu$ points. If, in addition, $\b\in
\nac c_\d$ for some $\d\in S^*$, then for every $x\in T(\b)$ define $y(x,\lng
P_{\d,\gamma}:\gamma\in \nac c_\d\wedge \gamma\le \b\rng)$ as follows:
let $\gamma_0=\max(c_\d\cap \b)$. When $\Lev(x)<\gamma_0$ set
$x_0$ as the supremum (in $T(\a)$) of $ \lng y(x,\lng P_{\d,\a}:\a\le
\a^*)\rng):\a^*\le\gamma_0\wedge \a^*\in\nacc c_\d\rng$; else,
$x_0=x$. As $|P_{\d,\b}|\le\cf\b$, we can in $\cf\b$ steps choose  a
cofinal branch above $x_0$ which has a point above an element from $A$
for every $A\in P_{\d,\b}$ which is a maximal antichain of $T_\b$. Let
the required $y$ be the supremum of this branch.

If $\b$ is a limit and $\cf\b=\mu$, distinguish two cases: case (a):
$\b=\d\in S^*$. So we should satisfy demand (vii), namely add bounds
precisely to those branches which for some $\d^*\ge\d$ in $S^*$, $\d\in\acc
c_{\d^*}$, are of the form $\lng y(x,t_\gamma):\gamma\in \nac
c_{\d^*}\cap \b\rng$ where  $t_\gamma=\lng
P_{\d^*,\zeta}:\zeta\le\gamma\wedge
\zeta\in \nac c_{\d^*}\rng$. By (f) this costs only the addition of
$\mu$ new elements. If, in addition, there is some $\d'\in S^*$
such that $\d\in \nac c_{\d'}$, we should define $y(x,\lng
P_{\d',\gamma}:\gamma\in \nac c_{\d'}\wedge \gamma\le\d\rng)$ for all
$x\in T(\d)$. This presents no problem: as $|P_{\d',\d}|<\mu$, we
attach to each $x$ some $x_0$ such that $x_0=x$ or $T_\d\models x_0>x$
and such that $x_0$ is above members from every maximal antichain in
$P_{\d',\d}$; now $y(x,\lng P_{\d',\gamma}:\gamma\in \nac
c_{\d'}\wedge \gamma\le\d\rng)$ will be the point in level $\d$ above
$x_0$ we obtained anyway to satisfy demand (vii).  

Case (b): $\cf\b=\mu$ and $\b\notin S^*$. Then when there is some $\d$
such that $\b\in\acc c_\d$ we realize enough limits to obtain
completeness under increasing sequences of the form $\lng
y(x,t_i):i<i^*\rng$.  By (f), we add thus $\le \mu$ elements. If
there is no such $\d$, just make sure, by adding $\mu$ points to the
tree in level $\b$, that above every $x\in T(\b)$ there is a point in
level $\b$.  This takes care also of (i). If there is some $\d'$ such
that $\b\in\nac c_{\d'}$, then for every $x\in T(\b)$ define $y(x,\lng
P_{\d',\gamma}:\gamma\in \nac c_{\d'}\wedge \gamma\le\d\rng)$ exactly
as in the case of smaller cofinality.\endproof 2

We will show now how to obtain from a non-reflecting stationary set a
special case of the prediction principle we used in the previous
theorem. One should substitute $P_{\d,\a}$ in the previous theorem by
$B_\a$ from the next theorem  to get the assumptions of the previous
theorem. 

\proclaim {3. Theorem}
Suppose $\l=\cf\l>\aleph_1$ , $S\su \l$ is stationary, non-reflecting,
and carries a diamond sequence $\lng A_\a:\a\in S\rng$, $S^*$ is a
given non reflecting stationary subset of $\l$, $S^*\cap S=\emptyset$
and $\d\in S^*\imply \cf\d>\aleph_0$. Then there are sequences $\ov C=\lng
c_\d : \d\in S^*\rng$ 
and $\ov B=\lng B_\a:\a\in S\rng$ such that:
\item{(i)} $B_\a\su  \a$;
\item{(ii)} $\sup
c_\d=\d$ and $c_\d$ 
is a closed set of limit ordinals; 
\item{(iii)} if $\d,\d^*\in S^*$ and $\d\in \acc c_{\d^*}$, then there is some
$\a<\d$ such that $c_{\d^*}\cap (\a,\d)=c_{\d}\cap (\a,\d)$;
\item{(iv)} for every club $E\su
\l$ and set $X\su\l$ there are stationarily many $\d\in S^*$ such that
$\d=\sup\{\a\in \nac c_\d: \a\in S\cap E
\wedge X\cap \a=A_\a\}$.

\proclaim{Proof} 
We fix  some 1-1  pairing function $\lng -,-\rng$ from $\l\times \omega_0$ onto $\l$
and let $A^\a_n=\{\b<\a:\lng
\b,n\rng\in A\}$.  We may assume that for every countable sequence $\ov
X=\lng X_n:n<\om\rng$ of subsets of $\l$ there are stationarily many
$\a\in S$ such that for every $n$, $X_n\cap \a=A^\a_n$. Denote
by $S(\ov X)$, for a (finite or infinite) sequence of subsets of $\l$ the
stationary set $\{\a\in S:\bigwedge_{n} X_n\cap\a=A^\a_n\}$.

To every limit $\a<\l$ we attach a club of $\a$, $e_\a$, satisfying
$e_\a\cap S=e_\a\cap S^*=\emptyset$, $\otp e_\a=\cf\a$ and $e_\a$
contains only limit ordinals whenever $\a\in S^*$. Let $\ov C_0=\lng
e_\d:\d\in S^*\rng$. Suppose that $\ov C_n=\lng c^n_\d:\d\in S\rng$ is
a bad candidate for the job, namely that there are a club $E_n$ and a
set $X_n$ such that for every $\d\in E_n\cap S^*$ the set
$\{\a\in\nacc c^n_\d:\a\in S(X_n)\cap E_n\}$ is bounded below $\d$.
(Surely, we may assume that $E_n$ is as thin as we like --- in
particular that all its members are limits).  Define $\ov
C_{n+1}$ by induction on $\d$: For every $\gamma\in c^n_\d$, we define 
$c^{n+1}_\d\cap (\gamma,\min c^n_\d\sm (\gamma+1))$ (where $(\gamma,\b)$
denotes, as usual, an open interval of ordinals), and we let
$c^{n+1}_\d=c^n_\d\cup\bigcup\{c^{n+1}_\d\cap(\gamma,\min
c^n_\d\sm(\gamma+1)):\gamma\in c^n_\d\}$. This is well defined, as every
$\gamma\in c^n_\d$ has a successor in $c^n_\d$. So denote by $\b$ the
ordinal $\min c^n_\d\sm(\gamma+1)$, and let

$$c^{n+1}_\d\cap (\gamma,\b) =
\cases{ 
c^{n+1}_\b\cap(\gamma,\b)& if $\b\in S^*$\cr
 & \cr
\emptyset&  if $\b\in S(X_0,\cdots,X_n)$\cr
 & \cr
\{\a:\gamma<\a<\b\wedge (\exists \zeta\in e_\b)(\a=\sup(\zeta\cap E_n))\}&
otherwise}  \leqno{ (*)}
$$

Note that for the definition to be consistent, $\b\in c^n_\d$ must
always be limit (and this is indeed the case). 

\proclaim {3.1 Lemma} Suppose that $\ov C_n$ is defined for $n\le m$.
Then for For every $n<m$ and $\d\in S^*$: 
\ss (0) If $\a\in c_\d^n$ then $\b$ is a limit ordinal;
\ss (1) $c^n_\d$ is closed.
\ss (2) $c^n_\d\su c^{n+1}_\d$.
\ss (3) If $\a\in S^*\cap \acc c_\d^{n}$, then $c^{n}_{\a}$
and $c^{n+1}_\d\cap\a$ have a common end segment.
\ss (4) If $\a\in c^{n+1}_\d\cap S(X_0,\cdots,X_n)$, then
$\a\in\nacc c^{n+1}_\d$.

\proclaim{Proof}: (2) is true by the definition of $c^{n+1}_\d$ for
every $n$ and $\d\in S^*$.  (0), (1), (3) and (4) are proved by induction
on $n$ and $\d$. 

For $n=0$ we know that $e_\d=c^0_\d$ is all limits
and is closed, so (0) and (1) hold. (3) is vacuously true, because
$e_\d\cap S^*=\emptyset$, and (4) is vacuously true because $e_\d\cap
S=\emptyset$. 

For  $n+1$: 

 (0): Suppose $\a\in c^{n+1}_\d$. If $\a\in c^n_\d$ then it is a limit
ordinal by (0) and the induction hypothesis on $n$. If $\a\notin
c^n_\d$, let $\gamma=\sup c^n_\d\cap \a$.  Because of (1) and the
induction hypothesis $\gamma<\a$. Let $\b=\min c^n_\d\sm(\a+1)$. If
$\b\in S^*$ then $c^{n+1}_\d\cap (\gamma,\b)=
c^{n+1}_\b\cap(\gamma,\b)$. So $\a\in c^{n+1}_\b$, and by the
induction hypotheses on $\b$, $\a$ is a limit ordinal. If $\b\notin
S^*$, the $c^{n+1}_\d\cap(\gamma,\b)=\{\a:\gamma<\a<\b,\quad (\exists
\zeta\in e_\b)(\a=\sup \zeta\cap E_n)\}$. Therefore for some $\zeta\in
e_\b$ our $\a$ is $\sup (\zeta\cap E_n)$. Since $E_n$ is a club,
$\a\in E_n$. But $E_n$ is a club of limits, so $\a$ is limit.

(1) Suppose that $\a\in \acc c^{n+1}_\d$, and we wish to show $\a\in
c^{n+1}_\d$. If $\a\in \acc c^n_\d$, then because of (1) and the
induction hypothesis on $n$ $\a\in c^n_\d$ and (by (2)) $\a\in
c^{n+1}_\d$. Else, $\gamma=\sup\a\cap c^n_\d$ and $\b=\min c^n_\d\sm
(a+1)$, $\gamma<\a<\b$. If $\b\in S^*$ then $\a\in \acc c^{n+1}_\b$. By
the induction hypothesis on $\b$ and (1), $\a\in c^{n+1}_\d$.
Otherwise, $\a$ is a limit of $\lng a_i:i<i^*\rng$ such that
$\a_i=\sup \zeta_i\cap E_n\in c^{n+1}_\d$. So clearly $\a\in E_n$. Let
$\zeta^*$ be the minimal in $e_\b$ above $\a$. So
$\a=\sup\zeta^*\cap E_n$. Therefore $\a\in c^{n+1}_\d$. 

Before	 proving (3) we note:

\proclaim{3.2 Fact} Suppose $\gamma\in c^n_\d$ and $\b=\min
c^n_\d\sm(\gamma+1)$. If $\b\notin S^*$ and  $\a\in
c^{n+1}_\d\cap(\gamma,\b)$ is a limit of $c^{n+1}_\d$, then $\a\in
e_\b$.

Indeed, if $\a=\sup\{\a(i):i<i^*\}$, where $\a(i)=\sup\zeta(i)\cap
E_n$ are elements in $c^{n+1}_\d$,  $\a\in E_n$. Therefore every
$\zeta(i)<\a$ (or else $\sup\zeta(i)\cap E_n\ge \a>\a(i)$). But
$\zeta(i)\ge \a(i)$, so $\a$ is a limit of $e_\b$. As $\a<\b$ and
$e_\b$ is closed, $\a\in
e_\b$.\endproof {3.2}

(3): Let $\a\in \acc c^{n+1}_\d\cap S^*$, and we wish to show that
$c^{n+1}_\d$ and $c^{n+1}_\a$ have a common end segment. If $\a\in \acc
c^n_\d$, then by the induction hypothesis on $n$ and (3), we know that
$c^n_\d$ and $c^n_\a$ have a common end segment; say they agree on
the interval $(\a(0),\a)$. This means in particular that for every
$\gamma\in c^n_\d\cap (\a(0),\a)$, $\a\in c^n_\a$ and $\min c^n_\d\sm
(\gamma+1) =\min c^n_\a\sm (\gamma+1)=:\b$. Therefore also $c_\d^{n+1}\cap
(\gamma,\b)=c^{n+1}_\a\cap(\gamma,\b)$, and consequently
$c^{n+1}_\d\cap(\a(0),\a)=c^{n+1}_\a\cap (\a(0),\a)$. So assume that
$\a\notin \acc c^n_\d$. The first possibility is that $\a\notin
c^n_\d$ altogether. In this case  let $\gamma<\a<\b$ assume their traditional
roles as the last ordinal of $c^n_\d$ below $\a$ and the first
above. If $\b\in S^*$, then by the induction hypothesis on $\b$ we
know that $c^{n+1}_\b$ and $c^{n+1}_\a$ have a common end segment; but
$c^{n+1}_\d\cap(\gamma,\b)=c^{n+1}_\b\cap(\gamma,\b)$, so it follows
that $c^{n+1}_\d$ and $c^{n+1}_\a$ have a common end segment. 

If
$\b\notin S^*$, then by the Fact  above, $\a\in e_\b$ ---
contradiction to $e_\b\cap S^*$ is empty. So this subcase is non
existent. 

The last case is: $\a\notin\acc c^n_\d$ but $\a\in c^n_\d$, or in short
$\a\in \nac c^n_\d$. Let $\gamma$ be the last element of $c^n_\d\cap\a$.
Then by $(*)$, $c^{n+1}_\d\cap
(\gamma,\a)=c^{n+1}_\a\cap(\gamma,\a)$. 

(4): Suppose that $\a\in S(X_0,\cdot,X_n)\cap
c^{n+1}_\d $.  We should see that $\a\in \nac c^{n+1}_\d$. Let $m\le
n+1$ be the minimal such that $\a\in c^m_\d$. It is enough to prove
that $\a\in \nac c^m_\d$, because by $(*)$ it is clear that if
$\a\in S(X_0,\cdots,X_m)\cap \nac c^m_\d$ then $\a$ will remain a
non-accumulation point in $c^{m+1}_\d$ (because nothing will be added in the
interval below it). So without loss of generality we may assume that
$\a\in c^{n+1}_\d\sm c^n_\d$. So denote by $(\gamma,\b)$, as usual,
the unique minimal interval with end points in $c^n_\d$ to which $\a$
belongs. First case: $\b\in S^*$. So $\a\in c^{n+1}_\b$; and by the
induction hypothesis on $\b$, $\a\in \nac c^{n+1}_\b$. So this case is
done. Otherwise, $\b\notin S^*$. So by the Fact above, if $\a$ were a
limit of $c^{n+1}_\d$, it would be in $e_\b$. But $\a\in S$, and therefore
cannot be in $e_\b$ by the very choice of $e_\b$. Therefore $\a\in \nac
c^{n+1}_\d$. (This is where the non reflection of $S$ is used in an
essential way).\endproof {3.1}

\proclaim{3.3 Claim}: There is some $n<\om$ for which $\ov C_n$ and
$\lng A^\a_n:\a\in S\rng$ are as
required.

\proclaim {Proof}: Suppose not. Let $\ov X_\om=\lng X_n:n<\om\rng$.
Let $E=\cap_n E_n$ and $E'=\acc (S(\ov X_\om)\cap E)$. So $E'$ is a
club. Pick some $\d\in S^*\cap E'$. For every $n$ there is a bound
below $\d$ of the set $\{\a\in \nacc c^n_\d:\a\in S(X_n)\cap E_n\}$.
As $\cf\d>\aleph_0$, let $\a^*<\d$ bound $\a(n)$ for all $n$. Let
$\d>\b>\a^*$ be in $S(\ov X_\om)\cap E$. So for every $n$,
$X_n\cap\b=A^\b_n$ and $\b\in E_n$. If $\b\in c^n_\d$ for some $n$,
then by (4)  $\b\in\nac c^n_\d$ --- a
contradiction to $\b>\a(n)$. So $\b\notin c^n_\d$ for all $n$.
Therefore for
every $n$ we may define $(\gamma(n),\b(n))$ as the minimal interval
with ends in $c^n_\d$ which contains $\b$.

\proclaim{3.4 Claim} $\b(n+1)<\b(n)$.

\proclaim{Proof} By its definition, $\b(n)\in\nac c^n_\d$. In the case
$\b(n)=\d^*\in S^*$, there are clearly elements in $c^{n+1}_{\d^*}$
above $\b$ and below $\b(n)$, so the claim is obvious. The case
$\b(n)\in S(X_0,\cdots,X_n)$ is impossible because of (4). In the
remaining case, $c^{n+1}_\d\cap
(\gamma(n),\b(n))=\{\a:\gamma(n)<\a<\b(n), (\exists \zeta\in
e_{\b(n)})(\a=\sup\zeta\cap E_n)\}$. Let $\zeta^*>\b$ be in
$e_{\b(n)}$. As $\b\in E\su E_n$, $\sup \zeta^*\cap E_n \ge \b$. But the
right hand side of this inequality belongs to $c^{n+1}_\d$, while
$\b$ does not; therefore $\sup \zeta^*\cap E_n>\b$. So we see that
there are elements of $c^{n+1}_\d$ in $(\b,\b(n))$, therefore the
least of them, namely $\b(n+1)$ is smaller than $\b(n)$.\endproof {3.4}

This is clearly a contradiction.  We conclude that after finitely
many steps, $\ov C_{n+1}$ cannot be defined due to the lack of a
counterexample. This means that $\ov C_n$ and $\lng
B_\a:\a\in S\rng$ where $B_\a=A^n_\a$ satisfy (i), (ii) and (iv). By
(3) above, they satisfy (iii) as well. \endproof {3.3}

This shows that after finitely many corrections all the requirements
are satisfied, and our theorem is proved.\endproof {3}

\proclaim{4. Theorem} If the $e_\d$ we
pick in the proof of Theorem 3 
satisfy the additional condition that for every $\gamma<\l$ the set
$\{e_\d\cap \gamma:\a\in
\l\;{\rm is}\;{\rm limit}\}$ has cardinality smaller then $\l$, then the
resulting good $\ov 
C= \lng c_\d:\d\in S^*\rng$ satisfies that for every $\gamma<\l$,
$|\{c_\d\cap \gamma:\d\in S^*\}|<\l$

\proclaim{Proof} Let $\gamma<\l$ be given. We must show that
$|\{c^n_\d\cap\gamma:\d\in S^*\}|\le \mu$. Let $N\prec H(\chi,\in)$ for
some large enough $\chi$, $|N|<\l$, $\gamma\su N$, $\gamma\in N$,
$\{ e_\a\cap \gamma:\a<\l\;{\rm is}\;{\rm limit}\}\su N$,  $\lng
e_\a:\e<\l\;{\rm is\; limit }\rng\in N$, 
and $E_n,X_n\in N$ for every $n$.

We shall see that for every $n$ and $\d$, $c^n_\d\cap \gamma\in N$.
Since $|N|<\l$, this is enough. 

First we notice that if $\d< \gamma$ then $c_\d\in N$ and by
elementarity also $c^n_\d\in N$ for every $n$. Now we use induction on
$n$ to show that for every $\d>\gamma$, $c^n_\d\cap \gamma\in N$.
For $n=0$: if $\d>\gamma$ then $c^0_\d\cap\gamma=e_\d\cap\gamma\in N$
by the assumptions on $N$.
For $n+1$ we use induction on $\d$. Suppose, then, that for all
$\d'<\d$, $c^{n+1}_{\d'}\cap \gamma\in N$. 

We need the easy

\proclaim{4.1 Fact} If $(\a_0,\a_1)$ is a minimal interval of $c^n_\d\cap
(\gamma+1)$ then $c^{n+1}_\d\cap (\a_0,\a_1)\in N$.

\proof  By $(*)$ above, the definition of $c^{n+1}_\d\cap
(\a_0,\a_1)$ depends only on $e_{\a_1}$, $E_n$ and (if case there is
such) $c^{n+1}_{\a_1}$. All these objects are in $N$, so also
$c^{n+1}_\d\cap (\a_0,\a_1)\in N$.\endproof {4.1}

Denote $\gamma(\d)=\sup
c^n_\d\cap\gamma$. So $\gamma(\d)\le\gamma$. If $\gamma(\d)=\gamma$,
then $c^{n+1}_\d\cap\gamma=c^n_\d\cap\gamma\cup \bigcup_I
c^{n+1}_\d\cap I$ where $I$ runs over all minimal intervals of
$c^n_\d\cap (\gamma+1)$. So by the Fact above we are done. Else,
$\gamma(\d)<\gamma$. In this case define $\b(\d)=\min c^n_\d\sm
\gamma$. If $\b(\d)=\gamma$ then again we are done by the Fact. The
remaining case is $\gamma(\d)<\gamma<\b(\d)$. By the same Fact,
$c^{n+1}_\d\cap \gamma(\d)\in N$. If $\b(\d)\in S^*$, then
$c^{n+1}_\d\cap (\gamma(\d),\gamma)=c^{n+1}_{\b(\gamma)} \cap
(\gamma(\d),\gamma)$. 
By the induction hypothesis, and since $\b(\d)<\d$, the latter set is
in $N$, and we are done. If $\b(\d)\notin S^*$, then either nothing is
added into $(\gamma(\d),\b(\d))$ (when $\b(\d)\in S(X_0,\cdots,X_n)$),
or $c^{n+1}_\d\cap
(\gamma(\d),\b(\d))=\{\a:\gamma(\d)<\a<\b(d)\,\;(\exists \zeta \in
e_{\b(\d)})(\a=\sup E_n\cap \zeta)\}$. So in this definition $N$ might not
know who $\b(\d)$ is, but $e_{\b(\d)}\cap \gamma\in N$. Therefore,
denoting by $\a^*$ the last member in $E_n\cap \gamma$, we can
determine in $N$ the set $c^{n+1}_\d\cap\a^*$. As to whether $\a^*$
itself is in this set or not, we need knowledge which is not available
in $N$, but who cares, as long as both possibilities are in $N$. \endproof 4

\proclaim {5. Corollary} If there is a non-reflecting stationary set
$S\su
\{\a<\mu^+:\cf\a<\mu\}$, and
$2^\mu=\mu^+$, $\mu^{<\mu}=\mu$, then there is a $\mu$-complete
Souslin tree on $\mu^+$.

\proclaim{6. Remark} This improves the result by Gregory in [G].

\proclaim{Proof} 
It is known (see [G] 2.1) that if $S\su \{\d\in\mu^+:\cf\d<\mu\}$ is
stationary, then $\mu=\mu^{<\mu}$ implies $\diamondsuit(S)$. As $S$ is
non reflecting, we can, for every limit $\a<\mu^+$, choose a closed
set $e_\a$, $\a=\sup e_\a$ and $\otp e_\a=\cf\a$ such that $e_\a\cap
S=\emptyset$. $\mu=\mu^{<\mu}$ implies that for every $\gamma<\mu^+$
the set $\{e_\a\cap\gamma:\a<\l,\;\a\;{\rm is \; limit}\}$ is of
cardinality at most $\mu$.  Use Theorem 3 and Theorem 4 to obtain the
assumptions of Theorem 2, $S$ being the given non reflecting
stationary set and $S^*$ being $\{\d<\l:\cf\d=\mu\}$. By Theorem 2 
there is an $\mu$-complete Souslin tree on $\mu^+$.\endproof 5

\proclaim {7. Problem} (1) Can the existence of such a tree be proved
in $ZFC+GCH$? 
\ss
(2) Can a Souslin tree on $\aleph_2$ be constructed from $GCH$ and
{\it two} stationary sets, each composed of ordinals of countable
cofinality, which do not reflect simultaneously? By
[Mg] this would raise the consistency strength of $GCH+SH(\aleph_2)$
to the consistency of the existence of a weakly compact cardinal.

\bigbreak
\centerline{\bf References}

[AbShSo] U.~Abraham, S.~Shelah and R.~M.~Solovay {\sl Squares with
diamonds and Souslin trees with special squares}, {\bf Fundamenta
Mathematicae} vol. {\bf 127} (1986) pp. 133--162

[De] K.~J.~Devlin, {\bf Aspects of Constructibility}, Lecture Notes
Math. {\it 354}, 240 pp., 1973. 

[DeJo] K.~J.~Devlin and H.~Johnsbr{\aa}ten {\bf The Souslin Problem},
Lecture Notes Math. {\it 405}, 132 pp., 1974.

[J] R.~Bj\"org Jensen, {\sl The fine structure of the constructible
hierarchy}, {\bf Annals of Mathematical Logic}, vol. {\bf 4} (1972), pp. 229--308.

[Je] T.~Jech, {\bf Set Theory}, Academic press, 1978.

[G] J.~Gregory, {\sl Higher Souslin Trees and the Generalized
Continuum Hypotheses}, {\bf The Journal of Symbolic Logic}, Volume
{\bf 41}, Number {\bf 3}, 1976, pp. 663--671.

[LvSh 104] R.~Laver and S.~Shelah, {\sl The $\aleph_2$-Souslin 
hypothesis}, 
{\bf Trans. of A.M.S.}, vol. {\bf 264} (1981), pp. 411--417.

[Mg] M.~Magidor, {\sl Reflecting stationary sets}, {\bf The Journal of
Symbolic Logic} vol {\bf 47}, number {\bf 4}, 1982 pp.755--771.

[Sh e] S.~Shelah, Non Structure Theory, {\bf Oxford University Press},
accepted. 

[Sh-g]  S.~Shelah, {\bf Cardinal Arithmetic}, submitted to {\bf
Oxford University Press}.

[Sh 365] S.~Shelah, {\sl There are Jonsson algebras in many 
inaccessible cardinals}, a chapter in {\bf Cardinal Arithmetic}. 

[ShSt 154]
S.~Shelah and L.~Stanley, {\sl Generalized Martin Axiom and the Souslin
Hypothesis for higher cardinality},
{\bf Israel J. of Math}, vol. {\bf 43} (1982) 225--236.

[ShSt 154a]
S.~Shelah and L.~Stanley, {\sl Corrigendum to "Generalized Martin's 
Axiom 
and  Souslin Hypothesis for Higher cardinality}, 
{\bf Israel J. of Math} vol. {\bf 53} (1986) 309--314.

[ShSt 279] S.~Shelah and L.~Stanley, {\sl Weakly compact cardinals
and non special Aronszajn trees} , {\bf Proc. A.M.S.}, vol. 104
no. {\bf 3} (1988) 887-897.
\end